\documentclass{siamart171218}
\usepackage{subfig}
\usepackage{forest}
\usepackage[section]{placeins}
\usepackage[export]{adjustbox}
\usepackage[title]{appendix}
\usepackage{amssymb,stmaryrd}
\usepackage{bm}

\usepackage{multirow}

\newcommand{\lp}{\left(}
\newcommand{\rp}{\right)}

\newcommand{\Z}{\mathbb{Z}}

\newcommand{\cat}[1]{\left \llbracket #1 \right \rrbracket}

\newcommand{\ghat}{\ensuremath{\hat{\bm{g}}}}
\newcommand{\fhat}{\ensuremath{\hat{\bm{f}}}}
\newcommand{\That}{\ensuremath{\hat{\bm{T}}}}
\newcommand{\Ahat}{\ensuremath{\hat{\bm{A}}}}
\newcommand{\ehat}{\ensuremath{\hat{\bm{e}}}}
\newcommand{\vhat}{\ensuremath{\hat{\bm{v}}}}
\newcommand{\uhat}{\ensuremath{\hat{\bm{u}}}}
\newcommand{\rhat}{\ensuremath{\hat{\bm{r}}}}

\newcommand{\Ren}{\text{Re}}

\newcommand{\pp}[2]{\frac{\partial #1}{\partial #2}}

\usepackage{lipsum}
\usepackage{amsfonts}
\usepackage{graphicx}
\usepackage{epstopdf}
\usepackage{algorithmic}
\ifpdf
  \DeclareGraphicsExtensions{.eps,.pdf,.png,.jpg}
\else
  \DeclareGraphicsExtensions{.eps}
\fi

\newsiamremark{remark}{Remark}
\newsiamremark{hypothesis}{Hypothesis}
\crefname{hypothesis}{Hypothesis}{Hypotheses}
\newsiamthm{claim}{Claim}

\headers{Schwarz Newton Krylov}{Kevin W. Aiton, Tobin A. Driscoll}

\title{Preconditioned nonlinear iterations for overlapping Chebyshev discretizations with independent grids \thanks{Submitted to the editors February 1, 2019.
\funding{This research was supported by National Science Foundation grant DMS-1412085.}}}

\author{Kevin W. Aiton, Tobin A. Driscoll}

\usepackage{amsopn}

\begin{document}

\maketitle

\begin{abstract}
  The additive Schwarz method is usually presented as a preconditioner for a PDE linearization based on overlapping subsets of nodes from a global discretization. It has previously been shown how to apply Schwarz preconditioning to a nonlinear problem. By first replacing the original global PDE with the Schwarz overlapping problem, the global discretization becomes a simple union of subdomain discretizations, and unknowns do not need to be shared. In this way restrictive-type updates can be avoided, and subdomains need to communicate only via interface interpolations. The resulting preconditioner can be applied linearly or nonlinearly. In the latter case nonlinear subdomain problems are solved independently in parallel, and the frequency and amount of interprocess communication can be greatly reduced compared to linearized preconditioning.   
\end{abstract}



\section{Introduction}
\label{sec:introduction}

Overlapping domain decomposition has been recognized as a valuable aid in solving partial differential equations since Schwarz first described his alternating method in 1870. (For straightforward introductions to the topic, see~\cite{Dolean2015,Smith2004}; for a more historical perspective, see~\cite{Gander2008}.) Overlapping decomposition provides a way to solve a problem on a global domain by exploiting its reduction to smaller subdomains. This creates geometric flexibility and allows special effort to be focused on small parts of the domain when appropriate. Domain decomposition also has a natural parallelism that is particularly attractive in the increasingly multicore context of scientific computing.

For a linear PDE, one typically seeks to apply a preconditioner for a Krylov iteration such as GMRES, in the form of solving problems on overlapping subdomains whose boundary data is in part determined by values of the solution in other subdomains. In the parallel context this is achieved by an additive Schwarz (AS) scheme. When one partitions the unknowns of a global discretization into overlapping subsets, the best form of AS are restricted AS (RAS) methods~\cite{Cai1999}, which do not allow multiple domains to update shared unknowns independently and thus over-correct. Typically, then, the subdomain problems are solved on overlapping sets, but the results are distributed in a nonoverlapping fashion. 

For nonlinear problems, the obvious extension of AS preconditioning is to apply it as described above on the linearized equations that are determined by a quasi-Newton iteration. We refer to this process as a \emph{Newton--Krylov--Schwarz} (NKS) procedure, reflecting the nesting order of the different elements of linearization, linear solver, and preconditioning.

Cai and Keyes~\cite{Cai2002} proposed instead modifying the \emph{nonlinear} problem using the Schwarz ansatz. In addition to yielding a preconditioned linearization for the Krylov solver, the preconditioned nonlinear problem exhibited more robust convergence for the Newton iteration than did the original nonlinear problem. They called their method ASPIN, short for \textit{additive Schwarz preconditioned inexact Newton}. As a technical matter, they did not recommend applying the true Jacobian of the system, preferring an approximation that required less effort. Subsequently, Dolean et al.~\cite{Dolean2016} pointed out that Cai and Keyes did not use the RAS form of AS preconditioning, and they proposed an improved variant called RASPEN that does. We refer to this type of nonlinear preconditioning as \emph{Schwarz--Newton--Krylov} (SNK), because the Schwarz ansatz is applied before the linearization begins.

Our interest is in applying the nonlinear preconditioning technique to spectral collocation discretizations in overlapping rectangles or cuboids, leading to globally smooth approximations constructed from a partition of unity~\cite{AitonTA}. In this context, there is not naturally a single global discretization whose degrees of freedom are partitioned into overlapping sets, because the Chebyshev (or Legendre, or other classical) nodes will not generally coincide within the overlapping regions. In principle one could link the degrees of freedom within overlap regions by interpolating between subdomains, but this process adds complication, computational time, and (in the parallel context) communication of data.

Here we present an alternative strategy that begins by replacing the original PDE problem with the Schwarz problems on the union of the subdomains. That is, rather than regarding the subdomains as solving the global PDE on a region that includes portions shared with other subdomains, each subdomain has a ``private copy'' of its entire region and is free to have its own solution values throughout. Of course, the new global problem is not solved until the interface values of every subdomain agree with values interpolated from other subdomains that contain the interface. As a Schwarz starting point, our technique has both NKS and SNK variants. 

One advantage of this new formulation is that interpolations need to be done only on lower-dimensional interfaces, rather than throughout the overlap regions. Another is that plain AS is preferred to RAS, because each subdomain has to update its own values separately. We show that it is straightforward to implement exact Jacobians for SNK with nothing more than the ability to do fully local PDE nonlinear and linearized solves, plus the ability to transfer values between subdomains through interface interpolations. We also derive a two-level method to prevent convergence degradation as the number of subdomains increases. The performance of the NKS and SNK methods is validated and compared through several numerical experiments.

\section{PDE problem and multidomain formulation}
\label{sec:pde}

The main goal of this work is to solve the PDE
\begin{subequations}
  \label{eq:global-problem}
  \begin{alignat}{2}
    \label{eq:pde}
    \phi(\bm{x},u) &= 0, & \qquad & \bm{x}\in\Omega, \\
    \label{eq:bc}
    \beta(\bm{x},u) &= 0, & \qquad & \bm{x}\in\partial\Omega,
  \end{alignat}  
\end{subequations}
where $u(\bm{x})$ is the unknown solution and $\phi$ and $\beta$ are nonlinear differential operators (with $\phi$ being of higher order). (We can easily extend to the case where $u$, $\phi$, and $\beta$ are vector-valued, but we use scalars to calm the notation.) Many Schwarz-based algorithms for~(\ref{eq:global-problem}) begin with a global discretization whose solution is accelerated by an overlapping domain decomposition. In this situation, some of the numerical degrees of freedom are shared by multiple subdomains---either directly or through interpolation---and proper use of additive Schwarz (AS) calls for the restricted-AS (RAS) implementation, which essentially insures that updates of shared values are done only once from the global perspective, not independently by the subdomains.

We take a different approach, replacing the original problem~(\ref{eq:global-problem}) with
\begin{subequations}
  \label{eq:overlapping-problem}
  \begin{alignat}{2}
    \label{eq:local-pde}
    \phi(\bm{x},u_i) &= 0, & \qquad & \bm{x}\in\Omega_i, \qquad i=1,\ldots,N,\\
    \label{eq:local-bc}
    \beta(\bm{x},u_i) &= 0, & \qquad & \bm{x}\in \Gamma_{i0}= \partial\Omega \cap \partial\Omega_i, \qquad i=1,\ldots,N,\\
    \label{eq:local-match}
    u_i & = u_j, & \qquad & \bm{x}\in\Gamma_{ij}= \partial\Omega_i \cap \Z_j, \qquad i,j=1,\ldots,N,
  \end{alignat}  
\end{subequations}
where now $u_1,\ldots,u_N$ are unknown functions on overlapping subdomains $\Omega_i$ that cover $\Omega$, and the $Z_i$ are nonoverlapping zones lying within the respective subdomains. Clearly any strong solution of~(\ref{eq:global-problem}) is also a solution of~(\ref{eq:overlapping-problem}), and while the converse is not necessarily so in principle,   we regard the possibility of finding a solution of~(\ref{eq:overlapping-problem}) that is not also a solution of~(\ref{eq:global-problem}) as remote in practice.

The key consequence of starting from~(\ref{eq:overlapping-problem}) as the global problem is that each overlapping region is covered separately by the involved subdomains; each is free to update its representation independently in order to converge to a solution. From one point of view, our discretizations of the overlap regions are redundant and somewhat wasteful. However, the fraction of redundant discrete unknowns is very modest. In return, we only need to interpolate on the interfaces, there is no need to use the RAS formulation, and the coarsening needed for a two-level variant is trivial (see section \ref{sec:NKS-two-level}).

\subsection{Discretization}

\label{sec:discretization}

We now describe a collocation discretization of~(\ref{eq:overlapping-problem}) for concreteness. Each subfunction $u_i(\bm{x})$ is discretized by a vector $\bm{u}_i$ of length $n_i$. By $\bm{u}=\cat{\bm{u}_i}$ we mean a concatenation of all the discrete unknowns over subdomains $i=1,\ldots,N$ into a single vector. Subdomain $\Omega_i$ is discretized by a node set $X_i\subset \overline{\Omega}_i$ and a boundary node set $B_i \subset \partial\Omega_i$. The total cardinality of $X_i$ and $B_i$ together is also $n_i$. The boundary nodes are subdivided into nonintersecting sets $G_{ij}=B_i\cap Z_j$ for $j\neq i$, and $G_{i0}= B_i \cap \partial\Omega$.

For each $i$, the vector $\bm{u}_i$ defines a function $\tilde{u}_i(\bm{x})$ on $\Omega_i$. These can be used to evaluate $\phi$ and $\beta$ from~(\ref{eq:overlapping-problem}) anywhere in $\Omega_i$. We define an $n_i$-dimensional vector function $\bm{f}_i$ as the concatenation of three vectors:
\begin{equation}
  \label{eq:fi}
  \bm{f}_i(\bm{u}_i) = 
  \begin{cases}
    \phi(\bm{x},\tilde{u}_i)& \text{ for all } \bm{x} \in X_i, \\
    \beta(\bm{x},\tilde{u}_i)& \text{ for all } \bm{x} \in G_{i0}, \\
    \tilde{u}_i(\bm{x}) & \text{ for all } \bm{x} \in G_{ij}, \; j=1,\ldots,i-1,i+1,\ldots,N.
  \end{cases}
\end{equation}
In addition, we have the linear \emph{transfer operator} $\bm{T}_i$ defined by
\begin{equation}
  \label{eq:Ti}
  \bm{T}_i\bm{u} = 
  \begin{cases}
    0 & \text{ for all } \bm{x} \in X_i, \\
    0 & \text{ for all } \bm{x} \in G_{i0}, \\
    \tilde{u}_j(\bm{x}) & \text{ for all } \bm{x} \in G_{ij}, \; j=1,\ldots,i-1,i+1,\ldots,N.
  \end{cases}
\end{equation}
Note that while $\bm{f}_i$ is purely local to subdomain $i$, the transfer operator $\bm{T}_i$ operates on the complete discretization $\bm{u}$, as it interpolates from ``foreign'' subdomains onto the parts of $B_i$ lying inside $\Omega$. 
Finally, we are able to express the complete discretization of~(\ref{eq:overlapping-problem}) through concatenations over the subdomains. Let $\bm{u}=\cat{\bm{u}_i}$, $\bm{f}(\bm{u})=\cat{\bm{f}_i(\bm{u}_i)}$, and $\bm{T}\bm{u} = \cat{\bm{T}_i\bm{u}}$. Then the discrete form of~(\ref{eq:overlapping-problem}) is the nonlinear equation
\begin{equation}
  \label{eq:discrete-full}
  \bm{f}(\bm{u}) - \bm{T} \bm{u} = \bm{0}.
\end{equation}
For a square discretization, the goal is to solve~(\ref{eq:discrete-full}), while in the least-squares case, the goal is to minimize $\bm{f}(\bm{u}) - \bm{T} \bm{u}$ in the (possibly weighted) 2-norm.

\subsection{Newton--Krylov--Schwarz}
\label{sec:NKS}

The standard approach to~(\ref{eq:discrete-full}) for a large discretization is to apply an inexact Newton iteration with a Krylov subspace solver for finding correcting steps from the linearization. Within the Krylov solver we have a natural setting for applying an AS preconditioner. Specifically, if we have a proposed approximate solution $\bm{u}$, we evaluate the nonlinear residual $\bm{r}=\bm{f}(\bm{u}) - \bm{T} \bm{u}$. We then (inexactly, perhaps) solve the linearization $\bigl[ \bm{f}'(\bm{u}) - \bm{T} \bm{u}\bigr] \bm{s}=-\bm{r}$ for the Newton correction $\bm{s}$, using a Krylov solver such as GMRES. These iterations are preconditioned by the block diagonal matrix $\bm{f}'(\bm{u})$, which is simply the block diagonal of the subdomain Jacobians $\bm{f}_i'(\bm{u}_i)$. We refer to this method as \emph{Newton--Krylov--Schwarz}, or NKS. 

Implementation of NKS requires three major elements: the evaluations of $\bm{f}(\bm{u})$ and $\bm{T} \bm{u}$ for given $\bm{u}$, the application of the Jacobian $\bm{f}'(\bm{u})$ to a given vector $\bm{v}$, and the inversion of $\bm{f}'(\bm{u})$ for given data. All of the processes involving $\bm{f}$ are embarrassingly parallel and correspond to standard steps in solving the PDE on the local subdomains. Each application of the transfer operator $\bm{T}$, however, requires a communication from each subdomain to its overlapping neighbors, as outlined in Algorithm~\ref{alg:transfer}. This step occurs once in evaluating the nonlinear residual and in every GMRES iteration to apply the Jacobian. In a parallel code, the communication steps could be expected to be a major factor in the performance of the method. 

\begin{algorithm}
  \caption{Apply transfer operator, $\bm{T}\bm{u}$.}
  \label{alg:transfer}
  \begin{algorithmic}
    \STATE Interpret input $\bm{u}$ as concatenated $\cat{\bm{u}_i}$.
    \FOR{$j = 1,\ldots,N$ (in parallel)}
    \FOR{all neighboring subdomains $i$}
    \STATE Evaluate $\tilde{u}_j$ at $\bm{x}\in G_{ij}$.
    \ENDFOR
    \ENDFOR
   \end{algorithmic}
\end{algorithm}

\subsection{Two-level scheme}
\label{sec:NKS-two-level}

As is well known~\cite{Dolean2015}, AS schemes should incorporate a coarse solution step in order to maintain convergence rates as the number of subdomains increases. The methods described above depend on the subdomain discretization sizes $n_i$ of the collocation nodes and solution representation, respectively. Now suppose we decrease the discretization sizes to $\hat{n}_i$, and denote the corresponding discretizations of~(\ref{eq:discrete-full}) by $\fhat(\uhat) - \That\uhat = \bm{0}$.  We can define a restriction operator $\bm{R}$ that maps fine-scale vectors to their coarse counterparts. This operator is block diagonal, i.e., it can be applied independently within the subdomains. We can also construct a block diagonal prolongation operator $\bm{P}$ for mapping the solution representation from coarse to fine scales. 

We are then able to apply the standard Full Approximation Scheme (FAS) using the coarsened problem~\cite{Brandt2011}. Specifically, we solve the coarse problem 
\begin{equation}
  \label{eq:FAS-residual}
  \fhat( \ehat + \bm{R}\bm{u} ) - \That\ehat - \fhat(\bm{R}\bm{u}) + \bm{R} ( \bm{f}(\bm{u}) - \bm{T}\bm{u} ) = \bm{0} 
\end{equation}
for the coarse correction $\ehat$, and define $\bm{c}(\bm{u})=\bm{P}\ehat$ as the FAS corrector at the fine level. The procedure for calculating $\bm{c}$ is outlined in Algorithm~\ref{alg:coarse-eval}.

\begin{algorithm}[tbp]
  \caption{Evaluate FAS correction $\bm{c}(\bm{u})$.}
  \label{alg:coarse-eval}
  \begin{algorithmic}
    \STATE Apply Algorithm~\ref{alg:transfer} to compute $\bm{T}\bm{u}$.
    \STATE Compute (in parallel) $\uhat = \bm{R} \bm{u}$.
    \STATE Compute (in parallel)  $\rhat = \bm{R}( \bm{f}(\bm{u}) - \bm{T}\bm{u} ) - \fhat(\uhat)$. 
    \STATE Solve equation~(\ref{eq:FAS-residual}) for $\ehat$.
    \STATE Compute (in parallel) the prolongation $\bm{P} \ehat$. 
  \end{algorithmic}
\end{algorithm}

We also require the action of the Jacobian $\pp{\bm{c}}{\bm{u}}=\bm{P} \pp{\ehat}{\bm{u}}$ on a given vector $\bm{v}$. It is straightforward to derive from~(\ref{eq:FAS-residual}) that
\begin{equation}
  \label{eq:FAS-jac}
   \bigl[ \fhat'( \hat{\bm{e}} + \bm{R}\bm{u} ) -\That \bigr] \pp{\hat{\bm{e}}}{\bm{u}}  = 
    -\bigl( \fhat'(\hat{\bm{e}} + \bm{R}\bm{u}) - \fhat'(\bm{R} \bm{u}) \bigr)\bm{R}
    - \bm{R} \bigl( \bm{f}'(\bm{u}) - \bm{T} \bigr).
  \end{equation}
Note that the matrix $\fhat'( \hat{\bm{e}} + \bm{R}\bm{u} )$ should be available at no extra cost from the end of the Newton solution of~(\ref{eq:FAS-residual}). Algorithm~\ref{alg:FAS-jac} describes the corresponding algorithm for computing the application of $\bm{c}'(\bm{u})$ to any vector $\bm{v}$. Even though $\bm{c}'$ is of the size of the fine discretization, the computation requires only coarse-dimension dense linear algebra. 

\begin{algorithm}[tbp]
  \caption{Apply Jacobian $\bm{c}'(\bm{u})$ for the FAS corrector to a vector $\bm{v}$.}
  \label{alg:FAS-jac}
  \begin{algorithmic}
    \STATE Apply Algorithm~\ref{alg:coarse-eval} to compute $\ehat$, $\uhat$, and the final value of $\Ahat=\fhat'(\ehat + \uhat)$.
    \STATE Apply Algorithm~\ref{alg:transfer} to compute $\bm{T}\bm{v}$.
    \STATE Set (in parallel) $\rhat = \bm{R} \bigl( \bm{f}'(\bm{u}) \bm{v} - \bm{T}\bm{v} \bigr)$ and
    $\vhat = \bm{R}\bm{v}$.
    \STATE Set (in parallel) $\hat{\bm{b}} =  \Ahat\vhat - \fhat'(\uhat)\vhat + \rhat$.
    \STATE Solve the linear system $ ( \Ahat-\That) \hat{\bm{y}} =- \hat{\bm{b}}$ for $\hat{\bm{y}}$.
    \STATE Compute (in parallel) $\bm{P} \hat{\bm{y}}$.
  \end{algorithmic}
\end{algorithm}

Finally, we describe how to combine coarsening with the preconditioned fine scale into a two-level algorithm. If we were to alternate coarse and fine corrections in the classical fixed-point form, 
\begin{align*}
  \bm{u}^{\dagger} &= \bm{u} + \bm{c}( \bm{u} ), \\
  \bm{u}^{\text{new}} &= \bm{u}^{\dagger} + \bm{f}( \bm{u}^{\dagger} ) - \bm{T} \bm{u}^\dagger,
\end{align*}
then we are effectively seeking a root of
\begin{equation}
  \label{eq:2level-residual}
  \bm{h}(\bm{u}) := \bm{c}( \bm{u} ) + (\bm{f}-\bm{T})( \bm{u} + \bm{c}( \bm{u} ) ).
\end{equation}
Finally, the Jacobian of the combined map is straightforwardly
\begin{equation}
  \label{eq:2level-jac}
  \bm{h}'(\bm{u}) = \bm{c}'( \bm{u} ) + (\bm{f}'-\bm{T})( \bm{u} + \bm{c}( \bm{u} ) )\cdot(\bm{I} + \bm{c}'( \bm{u} )).
\end{equation}
Thus the action of $\bm{h'}$ on a vector can be calculated using the algorithms for $\bm{c}'$, $\bm{f}'$, and $\bm{T}$.

\section{Preconditioned nonlinear iterations}
\label{sec:preconditioning}

As shown in section~\ref{sec:NKS}, the inner Krylov iterations of the NKS method are governed by the preconditioned Jacobian $\bm{I} - [\bm{f}'(\bm{v})]^{-1}\bm{T}$. Following the  observation of Cai and Keyes~\cite{Cai2002}, we next derive a method that applies Krylov iterations to the same matrix, but arising as the natural result of preconditioning the \emph{nonlinear} problem. Specifically, we precondition~(\ref{eq:discrete-full}) by finding a root of the nonlinear operator
\begin{equation}
  \label{eq:SNK}
  \bm{g}(\bm{u}) := \bm{u} - \bm{f}^{-1}(\bm{T}\bm{u}).
\end{equation}
Evaluation of $\bm{g}$ is feasible because of the block diagonal (that is, fully subdomain-local) action of the nonlinear $\bm{f}$. Since we are therefore applying the Schwarz ansatz even before linearizing the problem, we refer to the resulting method as \emph{Schwarz--Newton--Krylov} (SNK). 

We have several motivations for a method based on~(\ref{eq:SNK}). First, one hopes that the nonlinear problem, being a (low-rank) perturbation of the identity operator, is somehow easier to solve by Newton stepping than the original form is. Second, the inversion of $\bm{f}$ means solving independent nonlinear problems on the subdomains with no communication, which well exploits parallelism. Finally, the same structure means that problems with relatively small highly active regions could be isolate the need to solve a nonlinear problem to that region, rather than having it be part of a fully coupled global nonlinear problem.   

An algorithm for evaluating $\bm{g}$ is given in Algorithm~\ref{alg:SNK-eval}. It requires one communication between subdomains to transfer interface data, followed by solving (in parallel if desired) the nonlinear subdomain problems $\bm{f}_i$ defined in~(\ref{eq:fi}). Note that the local problem in $\Omega_i$ is a discretization of the PDE with zero boundary data on the true boundary $\Gamma_{i0}$ and values transferred from the foreign subdomains on the interfaces.

\begin{algorithm}
  \caption{Evaluate SNK residual $\bm{g}(\bm{u})$.}
  \label{alg:SNK-eval}
  \begin{algorithmic}
    \STATE Apply Algorithm~\ref{alg:transfer} to compute $\bm{T}\bm{u}$.
    \FOR{$i = 1,\ldots,N$ (in parallel)}
    \STATE Solve $\bm{f}_i(\bm{u}_i-\bm{z}_i)=\bm{T}_i\bm{u}$ for $\bm{z}_i$.
    \ENDFOR
    \STATE Return $\cat{\bm{z}_i}$. 
  \end{algorithmic}
\end{algorithm}

Using the notation of Algorithm~\ref{alg:SNK-eval}, we have that
\[
  \bm{f}_i'(\bm{u}_i-\bm{z}_i) \left[I - \pp{\bm{z}_i}{\bm{u}} \right] = \bm{T}_i,
\]
which implies that applying $[\partial \bm{z}_i/\partial \bm{u}]$ to a vector requires a single linear solve on a subdomain, with a matrix that is presumably already available at the end of the local Newton iteration used to compute $\bm{g}$. The process for applying $\bm{g}'(\bm{u})$ to a vector is outlined in Algorithm~\ref{alg:SNK-jac}. 
\begin{algorithm}
  \caption{Apply Jacobian $\bm{g}'(\bm{u})$ to vector $\bm{v}$ for the SNK problem.}
  \label{alg:SNK-jac}
  \begin{algorithmic}
    \STATE Apply Algorithm~\ref{alg:transfer} to compute $\bm{T}\bm{v}$.
    \FOR{$i = 1,\ldots,N$ (in parallel)}
    \STATE Solve $[\bm{f}_i'(\bm{u}_i-\bm{z}_i)] \bm{y}_i = \bm{T}_i\bm{v}$ for $\bm{y}_i$.
    \ENDFOR
    \STATE Return $\cat{\bm{v}_i - \bm{y}_i}$. 
  \end{algorithmic}
\end{algorithm}

\subsection{Two-level scheme}
\label{sec:SNK-two-level}

The SNK method can be expected to require coarse correction steps to cope with a growing number of subdomains. An obvious approach to incorporating a coarse-grid correction step is to apply FAS directly, i.e., using the analog of~(\ref{eq:FAS-residual}) with fine $\bm{g}$ and coarse $\ghat$. However, doing so means inverting $\ghat$, which introduces another layer of iteration in the overall process.

We have found that it is simpler and successful to apply the FAS correction in the form of the original NKS method, as given in section~\ref{sec:NKS-two-level}. All we need to do is replace $\bm{f}(\bm{u})-\bm{T}\bm{u}$ and $\bm{f}'-\bm{T}$ in~(\ref{eq:2level-residual}) and~(\ref{eq:2level-jac})  by $\bm{g}$ and $\bm{g}'$, respectively.

\section{Numerical experiments}
\label{sec:experiments}

For all experiments we compared three methods: NKS, SNK, and the two-level SNK2. The local nonlinear problems for SNK, and the coarse global problems in SNK2, were solved using \texttt{fsolve} from the Optimization Toolbox. Each solver used an inexact Newton method as the outer iteration, continued until the residual was less than $10^{-10}$ relative to the initial residual. For the inner iterations we used MATLAB's \texttt{gmres} to solve for the Newton step $s_k$ such that
\begin{equation}
  \| F(x_k) + F'(x_k)s_k \| \leq \eta_k \| F(x_k) \|	
\end{equation}
where $\eta_0 = 10^{-4}$ and
\begin{equation}
  \label{eq:forcing}
  \eta_k = 10^{-4} \lp \frac{\| F(x_k) \|}{\| F(x_{k-1})\|} \rp^2.	
\end{equation}
Given certain conditions on $F(x)$, if the intial solution is close enough to the true solution then this set of tolerances will yield a sequence with near q-2 convergence \cite{Eisenstat1996}.

\subsection{Regularized driven cavity flow}
\label{sec:cavity}
The first example is a regularized form of the lid-driven cavity flow problem~\cite{hirsch2007numerical}, where we replace the boundary conditions with infinitely smooth ones. Using the velocity-vorticity formulation, in terms of the velocity $u,v$ and vorticity $\omega$ on $\Omega = [0,1]^2$ we have the nondimensionalized equations
\begin{equation}
  \begin{aligned}
    -\Delta u - \frac{\partial \omega}{\partial y} &=0, \\
    -\Delta v + \frac{\partial \omega}{\partial x} &=0, \\
    -\frac{1}{\Ren} \Delta \omega + u \frac{\partial \omega}{\partial x} + v \frac{\partial \omega}{\partial y} &=0,
  \end{aligned}
  \label{reg_cavity_flow}
\end{equation}
where $\Ren$ is the Reynolds number. On the boundary $\partial \Omega$ we apply
\begin{equation}
  \begin{aligned}
    u &= \begin{cases}
      \exp \left(\frac{\left(\frac{y-1}{0.1}\right)^2}{1-\left(\frac{y-1}{0.1}\right)^2}\right), & y>0.9, \\
      0, & y \leq 0.9,
    \end{cases} \\
    v &= 0, \\
    \omega &= - \frac{\partial u}{\partial y} + \frac{\partial v}{\partial x},
  \end{aligned}
  \label{cavityflow_initalguess}
\end{equation}
similar to the boundary conditions in~\cite{shen1990numerical}.

We divided $\Omega$ into 16 overlapping patches of equal size (i.e. a 4 by 4 patch structure). Each subdomain was discretized by a second-kind Chebyshev grid of length 33 in each dimension. For the initial guess to the outer solver iterations we extended the boundary conditions (\ref{cavityflow_initalguess}) to all values of $x$ in the square $\Omega$. 

\begin{figure}
  \centering
  \subfloat[Nonlinear residuals with $Re = 100$]{
    \includegraphics[scale = 0.4]{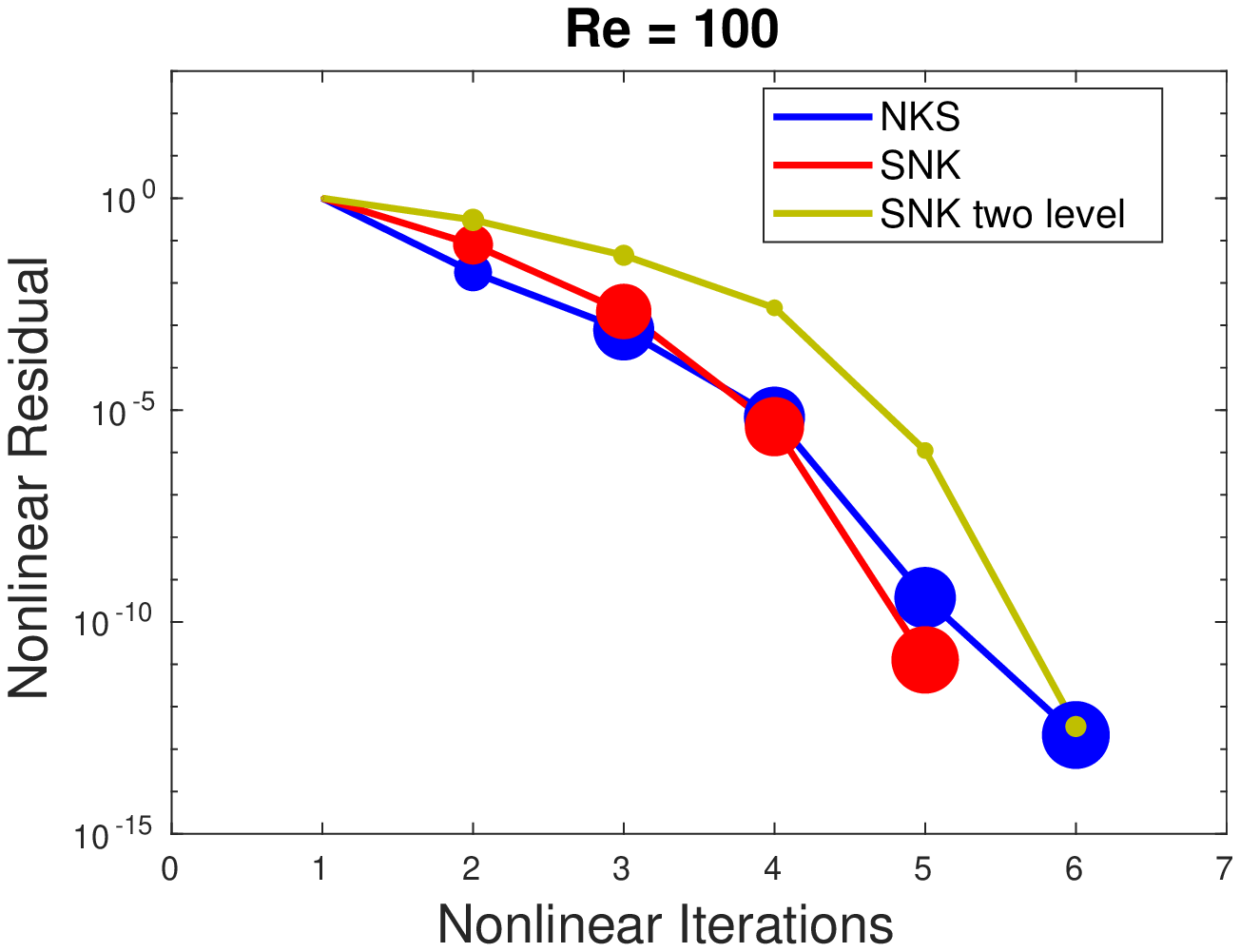}
    \label{cavityflow_100}
  }
  \subfloat[Nonlinear residuals with $Re = 1000$]{
    \includegraphics[scale = 0.4]{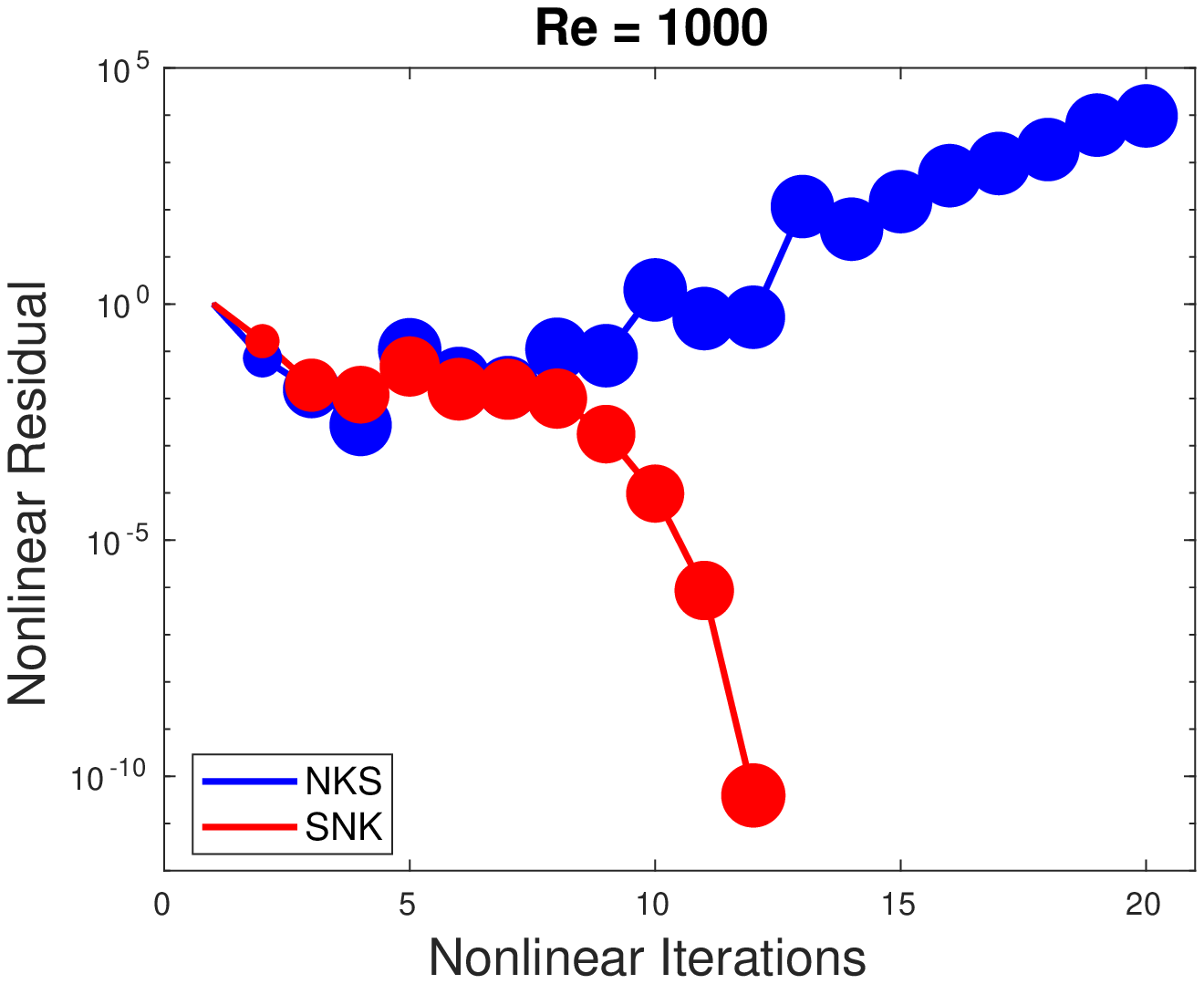}
    \label{cavityflow_1000}
  } 
  \caption{Nonlinear residuals, normalized by the residual of the initial guess, of the NKS, SNK, and SNK2 solvers on the regularized cavity flow problem (\ref{reg_cavity_flow})--(\ref{cavityflow_initalguess}). The area of each marker is proportional to the number of inner GMRES iterations taken to meet the inexact Newton criterion.}
  \label{reg_cavity_flow_residuals}
\end{figure}

The convergence of the three solvers is shown for $\Ren=100$, and SNK and NKS for $\Ren = 1000$ in Figure~\ref{reg_cavity_flow_residuals}. All three methods converge for $\Ren = 100$. The number of GMRES iterations per nonlinear iteration are similar for the NKS and SNK methods; this is to be expected since the linear system used to solve the Newton step is similar in both methods. We do however see a dramatic reduction in the number of GMRES iterations with the two-level SNK method. 

It is worth noting well that the computational time for each outer iteration varies greatly between solvers. The NKS residual requires only evaluating the discretized PDE and is thus is a good deal faster per iteration than the SNK solvers, which require solving the local nonlinear problems. In addition SNK2 must solve a global coarse problem in each outer iteration, but this added relatively little computing time. 

For the higher Reynolds number $\Ren = 1000$ we find that while SNK still converges, NKS does not, similar to what was reported in~\cite{Cai2002}. We also found that the coarse-level solver in SNK2 had trouble converging.

\subsection{Burgers equation}
\label{sec:burgers}

\begin{figure}
  \centering
  \includegraphics[scale = 0.6 ]{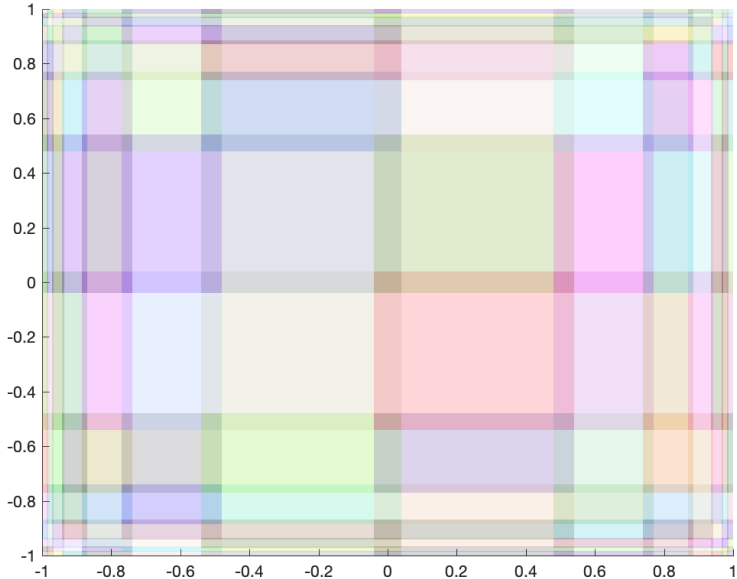}
  \caption{Subdomains for the Burgers experiments, found by adapting to the function  $\exp(-x^{20}/(1-x^{20}))exp(-y^{20}/(1-y^{20}))$ in order to increase resolution in the boundary layer.}
  \label{Burger_patch_structure}
\end{figure}

The second test problem is Burgers' equation,
\begin{equation}
  \nu \Delta u - u \cdot \nabla = 0,
  \label{simple_burgers}	
\end{equation}
on $\Omega = [-1,1]^2$, with Dirichlet boundary condition
\begin{equation}
  u = \arctan \lp \cos \lp \frac{3 \pi}{16} \rp x +\sin \lp \frac{3 \pi}{16} \rp y \rp.
  \label{burgers_boundary}
\end{equation}
This PDE was solved using a subdomain structure adapted to the function 
\[
  \exp\left( \frac{1}{1-x^{-20}} + \frac{1}{1-y^{-20}} \right)
\]
using the methods in \cite{AitonTA}, in order to help capture the boundary layers, as shown in Figure~\ref{Burger_patch_structure}. For the initial guess of the outer iterations, the boundary condition~(\ref{burgers_boundary}) was extended throughout $\Omega$.

\begin{figure}
  \centering
  \subfloat[Nonlinear residuals with $\nu = 1/400$]{
    \includegraphics[scale = 0.34]{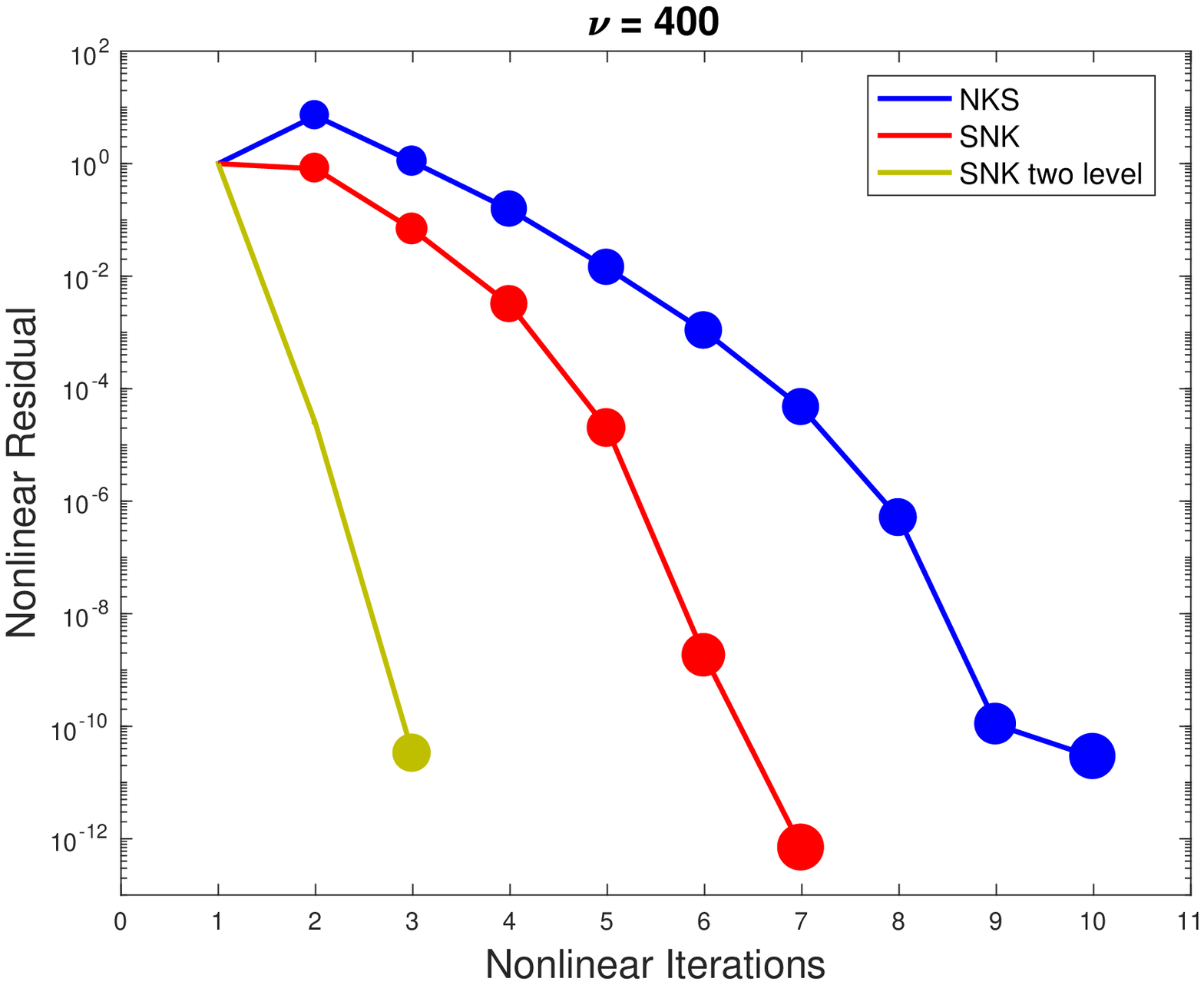}
    \label{burgers_200}
  }
  \subfloat[Nonlinear residuals with $\nu = 1/800$]{
    \includegraphics[scale = 0.34]{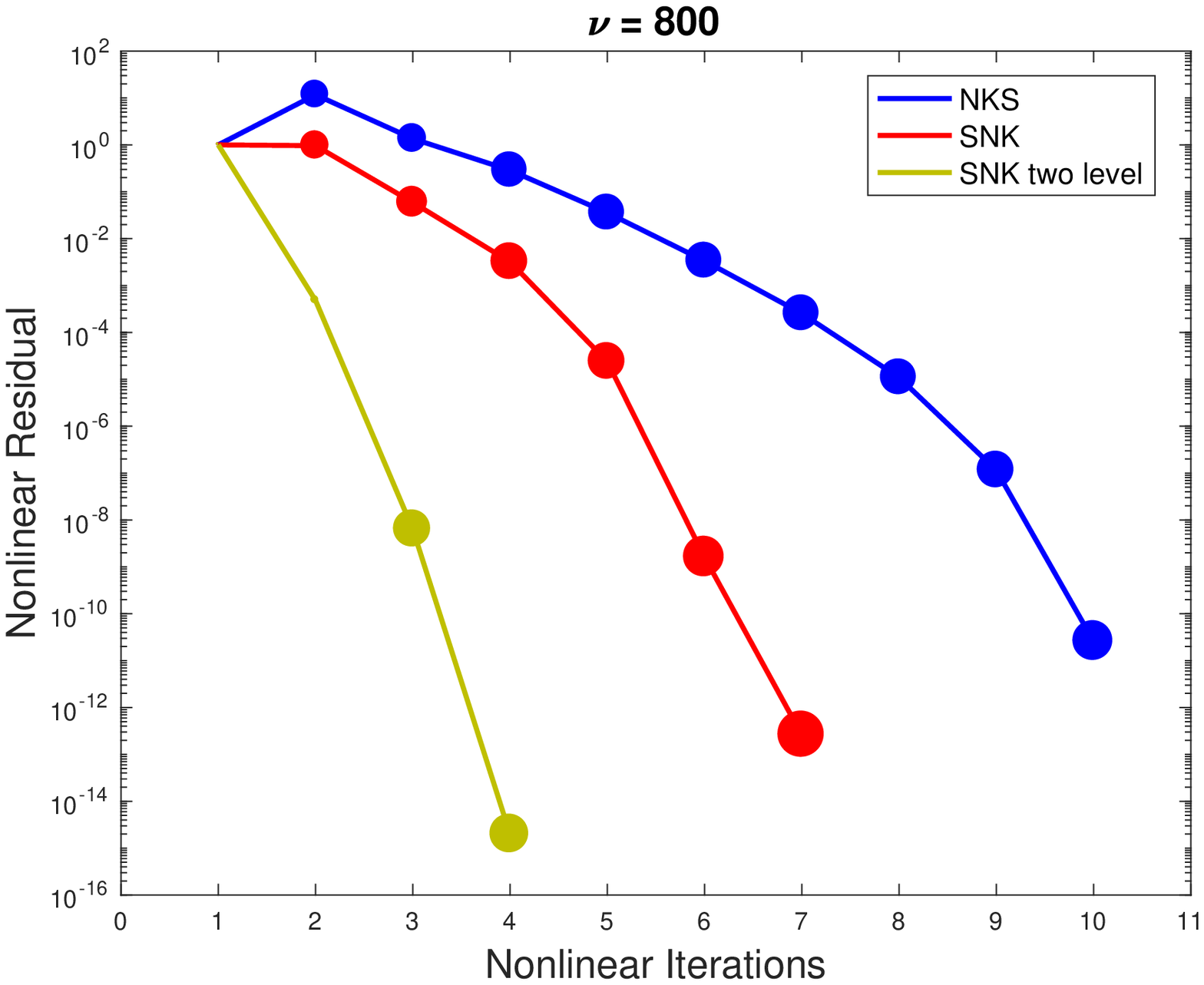}
    \label{burgers_800}
  }
  
    \subfloat[Nonlinear residuals with $\nu = 1/1000$]{
    \includegraphics[scale = 0.34]{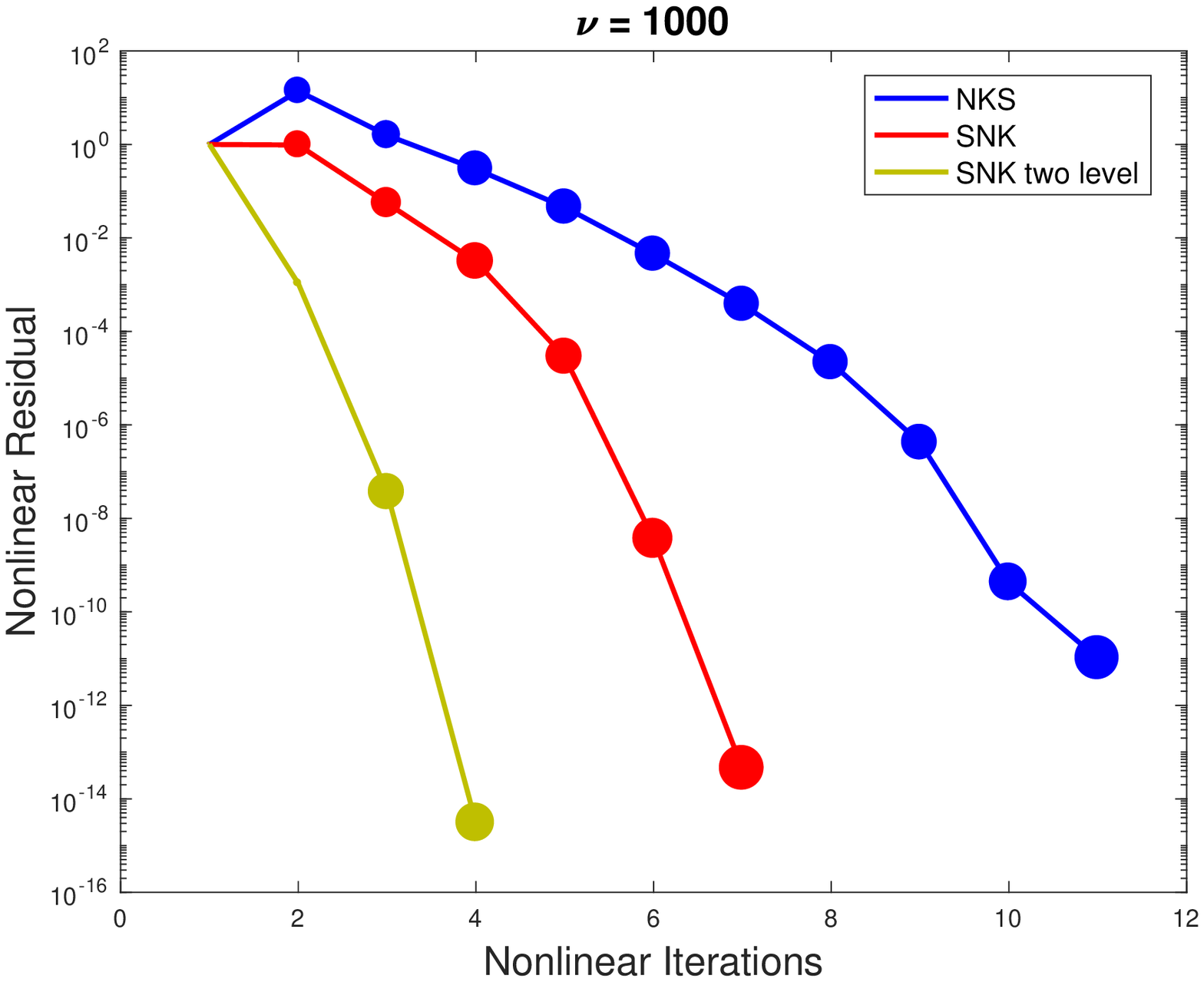}
    \label{burgers_1200}
  }
  \subfloat[Nonlinear residuals with $\nu = 1/1500$]{
    \includegraphics[scale = 0.34]{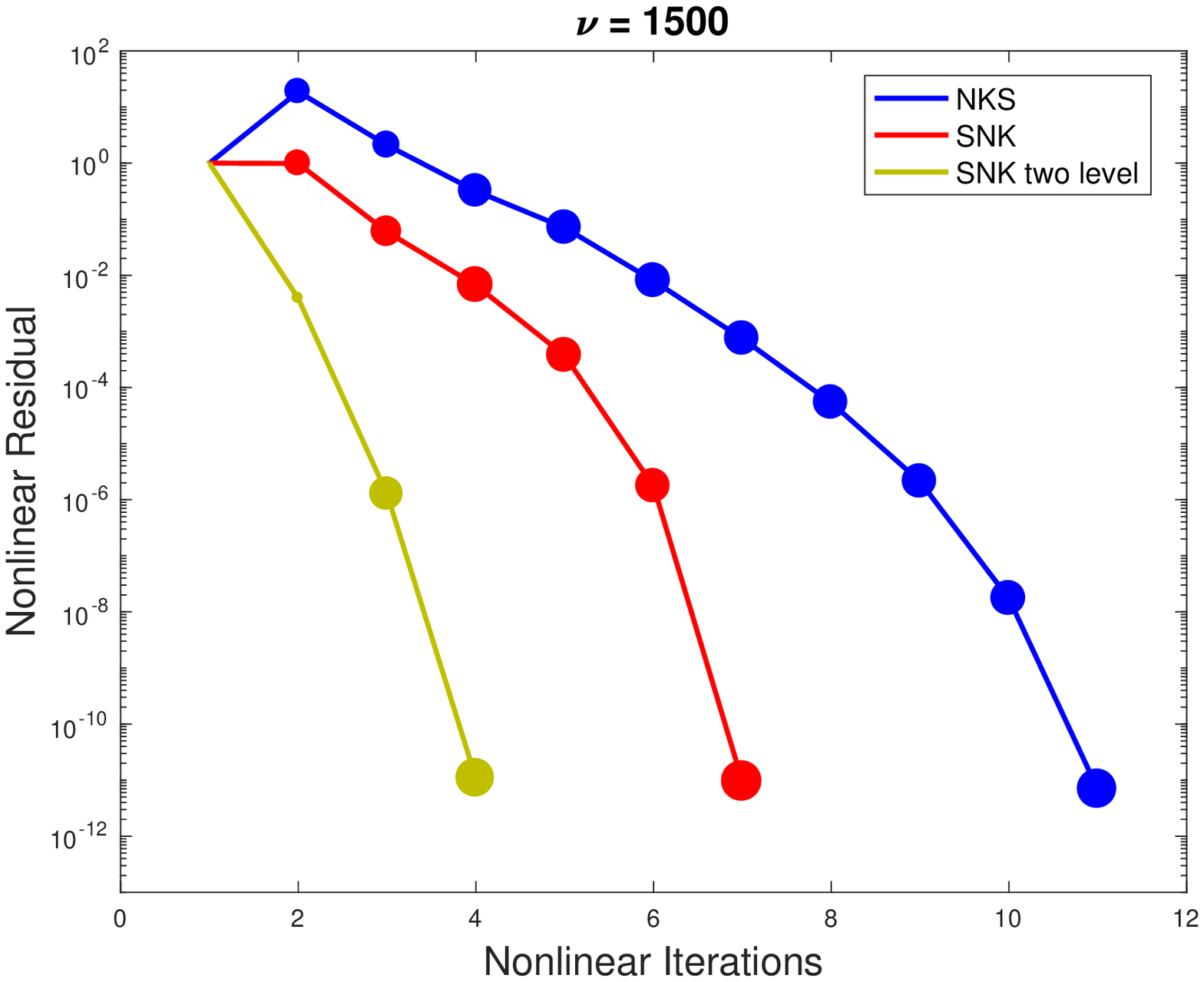}
    \label{burgers_1600}
  }
  \caption{Nonlinear residuals, normalized by the residual of the initial guess, of the NKS, SNK, and SNK2 methods to solve (\ref{simple_burgers})--(\ref{burgers_boundary}). The area of each marker is proportional to the number of inner GMRES iterations taken to meet the inexact Newton criterion.}
  \label{burgers_simp_test}
\end{figure}

Convergence histories for~(\ref{simple_burgers})--(\ref{burgers_boundary}) for $1/\nu = 400,800,1000,1500$ are given in Figure~\ref{burgers_simp_test}. We observe again that the SNK and SNK2 solvers seem quite insensitive to the diffusion strength, while the number of outer iterations in NKS increases mildly as diffusion wanes. Furthermore, SNK2 converges in about half as many outer iterations as SNK. 

\subsection{Parallel efficiency}
\label{sec:parallel}

A fully parallel-aware implementation of the methods would presumably distribute all the data and solving steps across cores, which would handle communication of interface values with neighbors when necessary. A simpler step was to modify our serial MATLAB implementation to use the \texttt{parfor} capability of the Parallel Computing Toolbox for the most compute-intensive loop in the SNK methods, that for the independent local nonlinear solves in the SNK residual evaluation, followed by factoring the final Jacobian matrices of these solutions. While we also tried parallelizing the loop for applying the inverses of the local Jacobians as part of the full Jacobian application in an inner Krylov iteration, the effect on timings was minimal or even detrimental due to the greater importance of communication relative to that computation.


The parallel SNK2 solver was applied to the Burgers experiment as described in section~\ref{sec:burgers}, but with the domain $\Omega$ split into an 8-by-8 array of uniformly sized subdomains, resulting in 3 outer iterations and a total of 18 inner iterations. The code was run on a compute node equipped with two 18C Intel E5-2695 v4 (for 36 total cores), 32 GB of DDR4 memory, and a 100 Gbps Intel OmniPath cluster network. The timing results for different numbers of parallel computing cores are given in Table~\ref{tab:parallel}. There is a good amount of speedup in evaluations of the nonlinear SNK residuals, consisting mainly of the solution of local nonlinear problems, which dominate the computing time for a small number of cores. However, the parallel efficiency is limited by the other parts of the implementation, most notably the Jacobian evaluations. 

\begin{table}
  \centering
  \begin{tabular}{rrrrrr}
    \multicolumn{1}{c}{Number} & \multicolumn{1}{c}{Total time} & \multicolumn{1}{c}{Speedup} & \multicolumn{1}{c}{Jacobian} & \multicolumn{1}{c}{Residual} & \multicolumn{1}{c}{Residual} \\
    \multicolumn{1}{c}{of cores} & \multicolumn{1}{c}{(sec.)} & \multicolumn{1}{c}{} & \multicolumn{1}{c}{time} & \multicolumn{1}{c}{time} & \multicolumn{1}{c}{speedup} \\
    \hline 
    1 & 62.3 & --- & 11.6 & 42.0 & --- \\
    2 & 43.5 & 1.43 & 12.7 & 24.1 & 1.74 \\
    4 & 33.4 & 1.86 & 11.3 & 15.5 & 2.71 \\
    6 & 30.0 & 2.08 & 11.5 & 12.1 & 3.47 \\
    8 & 28.3 & 2.21 & 11.3 & 10.5 & 4.00 \\
    12 & 27.6 & 2.26 & 11.9 & 9.3 & 4.51 \\
    16 & 28.6 & 2.18 & 13.2 & 8.8 & 4.78 \\
    20 & 28.2 & 2.21 & 13.1 & 8.3 & 5.04   
  \end{tabular}
  \caption{Parallel timing results for the Burgers equation experiment. ``Jacobian time'' is the total time spent within applications of the Jacobian to a given vector, and ``Residual time'' is the total amount of time spent evaluating the nonlinear SNK residual.}
  \label{tab:parallel}
\end{table}

\section{Discussion}
\label{sec:discussion}

We have described a framework for overlapping domain decomposition in which overlap regions are discretized independently by the local subdomains, even in the formulation of the global problem. Communication between subdomains occurs only via interpolation of values to interface points. This formulation makes it straightforward to apply high-order or spectral discretization methods in the subdomains and to adaptively refine them. 

The technique may be applied to precondition a linearized PDE, but it may also be used to precondition the nonlinear problem before linearization, to get what we call the Schwarz--Newton--Krylov (SNK) technique. In doing so, one gets the same benefit of faster Krylov inner convergence, but the resulting nonlinear problem is demonstrably easier to solve in terms of outer iterations and robustness. Although we have not given the derivation here, the Jacobian of the preconditioned nonlinear problem is readily shown to be a low-rank perturbation of the identity. Thus Kantorovich or other standard convergence theory for Newton's method~\cite{Dennis1987} may therefore suggest improved local convergence rates and larger basin of attraction. We have not yet pursued this analysis.

We have demonstrated that the SNK method can easily be part of a two-level Full Approximation Scheme in order to keep iteration counts from growing as the number of subdomains grows. The coarse level is simply a coarsening on each subdomain, so that restriction and prolongation steps can be done simply and in parallel. Indeed, the situation should make a fully multilevel implementation straightforward, as the multilevel coarsenings and refinements can all be done within subdomains.

The most time-consuming part of the SNK algorithm is expected to be typically in the solution of nonlinear PDE problems within each subdomain using given boundary data. These compute-intensive tasks require no communication and are therefore efficient to parallelize. By contrast, each inner Krylov iteration (i.e., Jacobian application) of both SNK and linearly preconditioned NKS requires a communication of interface data between overlapping subdomains, which appears to generate a more communication-bound form of parallelism. An additional feature of the SNK approach, mentioned also in~\cite{Cai2002}, is that subdomains of low solution activity can be expected to be found relatively quickly. We observed this to be the case in the cavity flow problem of section~\ref{sec:cavity}, where local solutions in regions of low activity were sometimes 3-4 times faster as those in regions with steep solution gradients. This presents a natural way to limit the spatial scope of difficult nonlinear problems, though it also raises questions for load balancing in a parallel environment.

Finally, we remark that an important extension in~\cite{AitonTA} is to use least-squares approximation rather than interpolation to incorporate nonrectangular (sub)domains. We have been able to write a least-squares (as opposed to collocation) generalization of SNK and test it in one dimension. We hope to make it the subject of future work.

\bibliographystyle{siamplain}
\bibliography{schwarz}

%
%
%

\end{document}